\documentclass[11pt]{amsart}

\usepackage{xypic}
\xyoption{all}
\usepackage{amstext}
\usepackage{amsthm}
\usepackage{amsopn}
\usepackage{amsfonts}
\usepackage{amsmath}
\usepackage{latexsym}
\usepackage{amscd}
\usepackage{amssymb}
\usepackage{amsmath}

\swapnumbers
\newtheorem{thm}[equation]{Theorem}
\newtheorem{prop}[equation]{Proposition}
\newtheorem{lem}[equation]{Lemma}
\newtheorem{cor}[equation]{Corollary}

\theoremstyle{definition}

\newtheorem{remark}[equation]{Remark}

\newtheorem{problem}[equation]{Problem}

\numberwithin{equation}{section}

\newcommand{\nsubsection}{
 \par\refstepcounter{equation}\medskip\noindent{\bf\theequation. }}

\def\limind{\mathop{\oalign{lim\cr
\hidewidth$\longrightarrow$\hidewidth\cr}}}
\def\limproj{\mathop{\oalign{lim\cr
\hidewidth$\longleftarrow$\hidewidth\cr}}}

\newcommand{\bbA}{{\mathbb A}}

\newcommand{\bbZ}{{\mathbb Z}}
\newcommand{\bbP}{{\mathbb P}}
\newcommand{\bbQ}{{\mathbb Q}}

\newcommand{\Spec}{\operatorname{Spec}}

\newcommand{\Hom}{\operatorname{Hom}}
\newcommand{\pr}{\operatorname{pr}}

\newcommand{\GL}{{\operatorname{GL}}}
\newcommand{\SL}{{\operatorname{SL}}}
\newcommand{\PSL}{{\operatorname{PSL}}}

\newcommand{\PGL}{{\operatorname{PGL}}}
\newcommand{\PGLn}{{\operatorname{PGL}_n}}

\newcommand{\Int}{{\operatorname{Int}}}
\newcommand{\oG}{\overline{G}}
\newcommand{\sdp}{\mathbin{{>}\!{\triangleleft}}}

\newcommand{\rank}{\operatorname{rank}}

\newcommand{\red}{\operatorname{red}}

\newcommand{\Char}{\operatorname{char}} 

\newcommand{\ed}{\operatorname{ed}}
\newcommand{\cd}{\operatorname{cd}}
\newcommand{\Gal}{\operatorname{Gal}}

\newcommand{\lra}{\longrightarrow}
\newcommand{\lga}{\longleftarrow}

\newcommand{\Br}{\operatorname{Br}}

\newcommand{\trdeg}{\operatorname{trdeg}}
\newcommand{\Sym}{{\operatorname{S}}}

\newcommand{\cal}{\mathcal}

\begin{document}
\title[$G$-torsors]{Resolving $G$-torsors by abelian base extensions}

\author{V. Chernousov}
\address{Department of Mathematics, University of Alberta,
    Edmonton, Alberta T6G 2G1, Canada}
\thanks{ V. Chernousov was supported by the Canada Research Chairs Program
and an NSERC Grant G121210944}
\email{chernous@math.ualberta.ca}

\author{Ph. Gille}
\address{Laboratoire de math\'ematiques, Universit\'e Paris-Sud, 91405 Orsay,
France}
\email{gille@math.u-psud.fr}

\author{Z. Reichstein}
\address{Department of Mathematics, University of British Columbia,
       Vancouver, BC V6T 1Z2, Canada}
\email{reichst\@@math.ubc.ca}
\urladdr{www.math.ubc.ca/$\stackrel{\sim}{\phantom{.}}$reichst}
\thanks{Z. Reichstein was partially supported by an NSERC research grant}

\subjclass{11E72, 14L30, 14E20}


\keywords{Linear algebraic group, torsor, non-abelian cohomology, 
$G$-cover, unramified cover, cohomological dimension, 
Brauer group, group action, fixed point obstruction}

\begin{abstract}
Let $G$ be a linear algebraic group defined over a field $k$.
We prove that, under mild assumptions on $k$ and $G$,
there exists a finite $k$-subgroup $S$ of $G$ such
that the natural map  $H^1(K, S) \lra H^1(K, G)$
is surjective for every field extension $K/k$.

We give several applications of this result 
in the case where $k$ an algebraically closed 
field of characteristic zero and $K/k$ is
finitely generated.  In particular, we prove
that for every $\alpha \in H^1(K, G)$
there exists an abelian field extension $L/K$
such that $\alpha_L \in H^1(L, G)$ is represented
by a $G$-torsor over a projective variety.
From this we deduce that $\alpha_L$ has trivial 
fixed point obstruction. We also show that
a (strong) variant
of the algebraic form of Hilbert's 13th problem
implies that the maximal abelian extension of $K$
has cohomological dimension $\le 1$. The last assertion,
if true, would prove conjectures of Bogomolov 
and K\"onigsmann, answer a question of Tits 
and establish an important case of Serre's 
Conjecture II for the group $E_8$.
\end{abstract}

\maketitle

\tableofcontents

\section{Introduction}

The starting point for this paper is the following
theorem, which will be proved in 
Sections~\ref{sect.reductive} and~\ref{sect.arb-field}.

\begin{thm} \label{thm4}
Let $G$ be a linear algebraic group defined over
a field $k$. Assume that one of the following conditions holds:

\smallskip
(a) $\Char(k) = 0$ and $k$ is algebraically closed, or

\smallskip
(b) $\Char(k) = 0$ and $G$ is connected, 

\smallskip
(c) $G$ is connected and reductive.

\smallskip
\noindent
Then there exists a finite $k$-subgroup $S$ of $G$, such that
the natural map $H^1(K, S) \lra H^1(K, G)$ is surjective
for every field extension $K/k$.
\end{thm}

Here, as usual $H^1(K, G)$ is the Galois cohomology set
$H^1(\Gal(\overline{K}/K), G)$; cf.~\cite{serre-gc}. Recall 
that this set does not, in general, have a group structure, 
but has a marked element, corresponding to the trivial
(or split) class, which is usually denoted by $1$.
Given a field extension $L/K$ we will, as usual,
denote the image of $\alpha$ under the natural
map $H^1(K, G) \lra H^1(L, G)$ by $\alpha_L$. 

In the course of the proof of Theorem~\ref{thm4} we 
will construct the finite group $S$ explicitly (see 
the beginning of Section~\ref{sect.reductive});
it is an extension of the Weyl group $W$ of $G$ by a finite
abelian group. (A related construction 
was used by Tits in~\cite{tits2}.)  Moreover, 
if $G$ is split and $k$ contains certain roots of unity 
then $S$ can be chosen to be a constant subgroup of $G$;
see Remark~\ref{rem.constant}.  
We also note that Theorem~\ref{thm4}(a) can be deduced from 
the results of Bogomolov (see~\cite[Lemma 7.3]{cs}); we are
grateful to J-L. Colliot-Th\'el\`ene for pointing this 
out to us.  We will include a self-contained proof 
of Theorem~\ref{thm4}(a) in Sections~\ref{sect.reductive}. 

In Section~\ref{sect.geom} we will discuss
Theorem~\ref{thm4}(a) in the context of invariant 
theory. In particular, we relate it to a result of 
Galitskii~\cite{galitskii} and use it to give a simple proof
of the no-name lemma, thus filling a small gap in the existing
literature; cf.~\cite[Section 4]{cs}.

Our other applications of Theorem~\ref{thm4} are motivated 
by the following question, implicit in the work of Tits~\cite{T}.

\begin{problem} \label{prob.tits} Let $G$ be a connected algebraic group
defined over an algebraically closed field of characteristic zero, 
$K/k$ be a field extension and $\alpha \in H^1(K, G)$. Is it true that
$\alpha$ can always be split by (i) a finite abelian field extension $L/K$
or (ii) by a finite solvable field extension $L/K$?
\end{problem}

Tits~\cite[Th\'eor\`eme 2]{T}
showed that Problem~\ref{prob.tits}(ii) has an affirmative answer
for every almost simple group of any type, other than $E_8$. 
(He also showed that for every such $G$, the solvable field 
extension $L/K$ can be chosen so that each prime factor of $[L:K]$ 
is a torsion prime of $G$.) Note that if Problem~\ref{prob.tits}(ii)
has an affirmative answer for fields $K$ of cohomological dimension $\le 2$, 
then we would be able to conclude, using an argument originally due to 
Chernousov, that $H^1(K, E_8) = \{ 1 \}$, thus proving an important 
(and currently open) case of Serre's Conjecture II; for details, 
see \cite[Chapter 6]{PR} or \cite[Th\'eor\`eme 11]{gille}.

We will say that $\alpha \in H^1(K, G)$ is {\em projective} 
if it is represented by a torsor over an irreducible complete variety $X/k$.  
In other words, $k(X) = K$, and $\alpha$ lies in the image 
of the natural map $H^1(X, G) \lra H^1(K, G)$, restricting
a torsor over $X$ to the generic point of $X$.
(Note that after birationally modifying $X$, 
we may assume it is smooth and projective.)
The split element of $H^1(K, G)$ is clearly projective, and
it is natural to think of projective elements of $H^1(K, G)$
as ``close" to being split.  The following result may thus be viewed 
as a ``first approximation" to the assertion of Problem~\ref{prob.tits}. 

\begin{thm} \label{thm2} Let
$k$ be an algebraically closed field of characteristic zero,
$G/k$ be a linear algebraic group, $K/k$ be a finitely generated
field extension, and $\alpha \in H^1(K, G)$.  
Then there exists a finite abelian extension $L/K$, such that
$\alpha_L$ is projective.
\end{thm}

Note that the group $G$ in Theorem~\ref{thm2}
is not assumed to be connected; in particular, the case where 
$G$ is finite (Proposition~\ref{prop.aby2}) is key
to our proof. On the other hand, in the case where $G$ is connected,   
 Theorem~\ref{thm2} does not imply an affirmative answer to
 Problem~\ref{prob.tits}. Indeed, while it is natural to think 
 of $\alpha_L$ as ``close to split", it may be not be literally 
 split, even in the case where $G$ is connected and simply connected.
 To illustrate this point, we will use a theorem of Gabber~\cite{cg}
 to construct a smooth projective 3-fold $X/k$  and a non-trivial 
 class $\alpha \in H^1(k(X), G_2)$ such that $\alpha$ is projective;
 see Proposition~\ref{prop.g2}. (Here $G_2$ denotes the (split) 
 exceptional group of type $G_2$ defined over $k$.) 

\smallskip
It is also natural to think of $\alpha \in H^1(K, G)$ as being 
``close to split" if $\alpha$ has {\em fixed point obstruction}; for
a precise definition, see Section~\ref{sect.fpo}. We will show 
that if $\alpha$ is projective then it has trivial fixed point 
obstruction; see Proposition~\ref{prop.fpo}. Combining this result with 
Proposition~\ref{prop.g2} yields another ``approximation" 
to the assertion of Problem~\ref{prob.tits}. 

\begin{cor} \label{cor.fpo}
Let $k$ be an algebraically closed field of characteristic zero,
$G/k$ be a linear algebraic group, $K/k$ be a finitely generated
field extension, and $\alpha \in H^1(K, G)$.  
Then there exists a finite abelian extension $L/K$, such that
$\alpha_L$ has trivial fixed point obstruction.
\qed
\end{cor}

In Section~\ref{sect.hilbert} we will use Theorem~\ref{thm4}(a) to relate
Problem~\ref{prob.tits} to a (strong) variant 
of Hilbert's 13th problem (Problem~\ref{prob.hilbert}).
We will show that if Problem~\ref{prob.hilbert} had an affirmative 
answer then so would Problem~\ref{prob.tits} (and, in fact, a much
stronger assertion would then hold; see Theorem~\ref{main} 
and Remark~\ref{rem.main}).

\section{Proof of Theorem~\ref{thm4}(a)}
\label{sect.reductive}

We begin with the following observation.
Let $k$ be a field of characteristic zero, $G/k$ be a linear 
algebraic group, and $R_u(G)$ be the unipotent radical of $G$.
Recall that $G$ has a Levi decomposition, $G = R_u(G) \sdp G_{\red}$,
where $G_{\red}$ is a reductive subgroup of $G$, uniquely determined
up to conjugacy. As usual, we shall refer to $G_{\red}$ as 
a {\em Levi subgroup} of $G$.


\begin{lem} \label{lem2b}
Let $i \colon G_{\red} \hookrightarrow G$ be a Levi subgroup of $G$. 
Then for any field extension $K/k$, the natural map
\[ i_* \colon H^1(K, G_{\red}) \lra H^1(K, G) \] 
is a bijection.
\end{lem}

\begin{proof} Let $\pi \colon G \lra G/R_u(G)$ be the natural projection.
By the Levi decomposition,
$G_{\red} \stackrel{i}{\hookrightarrow}
G \xrightarrow{\pi} G/R_u(G)$ 
is an isomorphism between $G_{\red}$ and $G/R_u(G)$.  Thus
\[ H^1(K, G_{\red}) \xrightarrow{i_*}
H^1(K, G) \xrightarrow{\pi_*} H^1(K, G/R_u(G)) \]
is a bijection between $H^1(K, G_{\red})$ and $H^1(K, G/R_u(G))$. 
By~\cite[Lemma 1.13]{sansuc}, $\pi_*$ is also a bijection. Hence, 
so is $i_*$.
\end{proof}

\begin{remark} \label{rem.red} 
Lemma~\ref{lem2b} tells us that if the natural map
\[ H^1(K, S) \lra H^1(K, G_{\red}) \]
is surjective then so is the natural map
\[ H^1(K, S) \lra H^1(K, G) \, . \]
In particular, in the course of proving Theorem~\ref{thm4}(a) and (b) 
we may replace $G$ by $G_{red}$ and thus assume that $G$ is reductive. 
\end{remark}

\def\limind{\mathop{\oalign{lim\cr
\hidewidth$\longrightarrow$\hidewidth\cr}}}
\def\limproj{\mathop{\oalign{lim\cr
\hidewidth$\longleftarrow$\hidewidth\cr}}}

We now proceed with the proof of Theorem~\ref{thm4}(a).
Let $k$ be an algebraically closed field of characteristic zero
and $G$ be a linear algebraic group defined over $k$. 
In view of Remark~\ref{rem.red}, we will assume that $G$ 
(or equivalently, the connected component $G^0$ of $G$) is reductive. 

Let $T$ be a maximal torus of $G$ and set $N=N_G(T)$ and $W=N_G(T)/T$.
Then $W$ is a finite group and $N$ is an extension of $W$ by $T$.
Let $\mu'={_nT}$ (resp. $\mu={_{n^2}T}$)  be the group 
of $n$-torsion points of $T$,
where $n = |W|$.  Consider the exact sequences
\[
1\to T \to N \stackrel{p}{\to} W \to 1
\quad \quad  \text{and} \quad \quad
1 \to  \mu' \to T \buildrel \times n \over \to T
\to 1
\, . \]
The first sequence yields a class in $H^2(W,T)$.
Since $n \cdot H^2(W,T)=0$, the second sequence tells us that
this class comes from $H^2(W, \mu')$. In terms of group extensions,
it means that there exists an extension $S'$
of $W$ by $\mu'$ such that $N$ is the push-out of
$S'$ by the morphism $\mu' \hookrightarrow T$.

In the same way, we obtain a group extension $S \subset N$ of $W$ by $\mu$.
Note that  $S'\subset S$ and $|S|= |W|^{\rank(G)^2 + 1}$.
Theorem~\ref{thm4}(a) is now an immediate consequence of the following
proposition. 

\begin{prop} \label{prop1}
Assume $G$ is reductive and $S$ is the finite subgroup of $G$
constructed above. Then the map $H^1(K, S) \rightarrow
H^1(K,G)$ is surjective for any field extension $K/k$.
\end{prop}

\begin{proof} We claim that the natural map
$H^1(K,N) \rightarrow H^1(K,G)$ is surjective for every
field extension $K/k$. Indeed, let
$\overline{K}$ be an algebraic closure of $K$.
For any $[z] \in H^1(K,G)$ the twisted group
$_zG^0$ is reductive and has a maximal torus $Q$.
Viewing $Q$ and $T$ as maximal tori in
$G^0(\overline{K})$, we see that they are
$\overline{K}$-conjugate; the claim now
follows from~\cite[Lemma III.2.2.1]{serre-gc}.

It remains to prove that the map $H^1(K, S)
\rightarrow H^1(K,N)$ is surjective. We will do this
fiberwise, with respect
to the map $p_*: H^1(K,N) \rightarrow H^1(K,W)$.
Fix $[b] \in H^1(K, N)$; our goal is to show that
$[b]$ lifts to $H^1(K, S)$.

\begin{lem} \label{lem2.4}  
Let $[a]=p_*([b]) \in H^1(K,W)$.
Then $[a]\in {\rm Im}( H^1(K,S) \buildrel q_* \over \to H^1(K,W))$.
\end{lem}

\begin{proof}
The obstruction to lifting $[a]$ to $H^1(K, S)$ is the class 
$$
\Delta([a]) \in  H^2(K, {_a\mu}),
$$
where ${_a\mu}$ denotes the group $\mu$, twisted
by the cocycle $a$ and $\Delta$ is the connecting map;
see ~\cite[I.5.6]{serre-gc}.
We now use the commutative diagram of horizontal 
exact sequences
\begin{equation} \label{diagr.lem2.4}
\begin{CD}
1 @>>> \mu' @>>> S' @>{q'}>> W @>>> 1 \\
&&  \cap &&  \cap && \mid \, \mid \\
1 @>>> \mu @>>> S @>{q}>> W @>>> 1 \\
&&  \cap &&  \cap && \mid \, \mid \\
1 @>>> T @>>> N @>{p}>> W @>>> 1 \\
\end{CD}
\end{equation}
and the functoriality of the connecting map $\Delta$. 
The obstruction to lifting
$[a]$ to $H^1(K, S')$, via $q'_*: H^1(K, S') \rightarrow H^1(K,W)$, 
is $\Delta'([a]) \in  H^2(K, _a\mu')$, where
$\Delta([a])$ is the image of $\Delta'([a])$ under  
the natural map $H^2(K, {_a\mu'}) \to H^2(K, {_a\mu})$.

The commutative exact diagram
\[ 
\begin{CD}
1 @>>> {_a\mu'} @>>> {_aT} @>{\times n }>> {_aT} @>>> 1 \\
&&  \cap && \mid \, \mid   && @V{\times n}VV \\
1 @>>> {_a\mu} @>>> {_aT} @>{\times n^2}>> {_aT} @>>> 1 \\
\end{CD}
\]
gives rises to the commutative exact diagram
$$
\begin{CD}
H^1(K, \, {_aT}) @>\Delta' >> H^2(K, \, _a\mu') @>>> H^2(K, \, {_aT})   \\
 @V{\times n}VV  @VVV      \mid \, \mid    \\
H^1(K, \, {_aT}) @>\Delta>> H^2(K, \, {_a\mu}) @>>> H^2(K, \, {_aT}) .
\end{CD}
$$
We shall now analyse this diagram. Recall that the middle vertical map
sends $\Delta'([a])$ to $\Delta([a])$. 
Since we are assuming that $[a]$ lifts to $[b] \in H^1(K, N)$, we have  
$$ \Delta'([a]) \in  \ker\Bigl( H^2(K, {_a\mu'})\to H^2(K, {_aT}) \Bigr).  $$
In other words, $\Delta'([a])$ lies in the image of $H^1(K, \, {_aT})$. 
Thus in order to prove the lemma (i.e., to prove that $\Delta([a]) = 0$),
it suffices to show that the left vertical map
$$
\begin{CD}
H^1(K, \, {_aT})   \\
 @V{\times n}VV    \\
H^1(K, \, {_aT}) 
\end{CD}
$$
is trivial.

Indeed, the torus $_aT$ is split by the Galois extension $L/K$ given by
$[a] \in H^1(K,W)= \Hom_{ct}( \Gal(\overline{K}/K), W)/ \Int(W)$;
the degree of this extension divides $n = |W|$.
The restriction-corestriction formula
\[ \times [L: K]= \operatorname{Cor}_K^L \circ 
\operatorname{Res}_K^L \]
and the fact that $H^1(L,T)=0$ (Hilbert's Theorem 90) imply that the map
\[ \times [L:K] \colon H^1(K,{_aT}) \rightarrow H^1(K,{_aT}) \]
is trivial.  Since $[L:K]$ divides $n$, the map
$\times n \colon  H^1(K,{_aT}) \rightarrow H^1(K,{_aT})$ is trivial
as well.
\end{proof}

We are now ready to finish the proof of Proposition~\ref{prop1}.
Let $[c] \in H^1(K,S)$ be such that $q_*([c])=[a]$.
The bottom two rows of~\eqref{diagr.lem2.4} give 
rise to the diagram 
\[ \begin{CD}
H^1(K, {_c \mu}) @>>> q_*^{-1}(a) \subset H^1(K, S) \\
 @VVV  @VVV  \\
H^1(K, {_c T}) @>>> p_*^{-1}(a) \subset H^1(K, N) . 
\end{CD} \]
Recall that our goal is to show that 
$[b] \in p_*^{-1}([a]) \subset H^1(K, N)$ lies 
in the image of $H^1(K, S)$.
A twisting argument~\cite[I.5.5]{serre-gc} shows that
the map $$H^1(K, {_{c}T}) \rightarrow p_*^{-1}([a])$$ 
is surjective; see~\cite[I.5.5]{serre-gc}.
Thus it suffices to prove that the vertical map
\[ \begin{CD}
H^1(K, {_{c} \mu}) \\
@VVV \\
H^1(K, {_{c}T}) \end{CD} 
\]
is surjective as well.  The cokernel of this map is given 
by the exact sequence
$$
H^1(K, {_{c}\mu}) \rightarrow  H^1(K,{_{c}T}) \buildrel
\times n^2 \over \rightarrow  H^1(K,{_{c}T}) .
$$
As we saw at the end of the proof of Lemma~\ref{lem2.4}, the map
$\times n : H^1(K,{_{c}T}) \rightarrow H^1(K,{_{c}T})$ is trivial
and hence, so is
$\times n^2 : H^1(K,{_{c}T}) \rightarrow H^1(K,{_{c}T})$.
We conclude that the map $H^1(K, {_{c}\mu}) \rightarrow  H^1(K,{_{c}T})$
is surjective. This completes the proof of Proposition~\ref{prop1}.
\end{proof}

\section{Proof of Theorem~\ref{thm4}(b) and (c)}
\label{sect.arb-field}

In view of Remark~\ref{rem.red} part (b) follows from part (c).
The rest of this section will be devoted to proving part (c).
We will consider three cases.

\smallskip
{\it Case $1$}. Let $G$ be a quasi-split adjoint group.
We denote by $T$ a maximal quasi-split torus in $G$,
$N=N_G(T)$ and $W=N_G(T)/T$. For every root $\alpha\in \Sigma=\Sigma(G,T)$,
where $\Sigma$ is the root system of $G$ with respect to $T$, the
corresponding subgroup $G_{\alpha}\le G$ is isomorphic (over a separable
closure of $k$) to either $\SL_2$
or $\PSL_2$.

Let $T_{\alpha}=T\cap G_{\alpha}$ and let
$w_{\alpha}\in N_{G_{\alpha}}(T_{\alpha})$ be a representative of
the Weyl group of $G_{\alpha}$ with respect to the maximal torus $T_{\alpha}$
given by a matrix
$$
\left(
\begin{array}{cc}
0 & 1\\
-1 & 0
\end{array}
\right).
$$
By Galois' criteria for rationality, the group $L$ generated by all
$w_{\alpha}$ is $k$-defined. One easily checks that the intersection
$L\cap T$ belongs to the $2$-torsion subgroup of $T$; in particular, 
$L$ is finite. 

Let $\mu={_{n^2} T}$ be
the $n^2$-torsion subgroup of $T$ where $n$ is the cardinality
of the Weyl group $W$. Consider the subgroup $S$ of $N$ generated
by $L$ and $\mu$. Now, arguing as in the proof of Proposition~\ref{prop1},
one checks that the canonical map $H^1(K,S)\to H^1(K,N)$
is surjective for every extension $K/k$. In the course
of the proof of Proposition~\ref{prop1} we showed that 
$H^1(K,N)\to H^1(K,G)$ is surjective. Then the composite map 
$H^1(K,S)\to H^1(K,N)\to H^1(K,G)$ is surjective as well.

\smallskip
{\it Case 2.} Let $G$ be an adjoint $k$-group. Denote by $G_0$
the quasi-split adjoint group of the same inner type as $G$. One knows
(see \cite{tits})
that $G={_a(G_0)}$ is the twisted form of $G_0$ for an appropriate
cocycle $a\in Z^1(k,G_0)$. If $S_0$ is the subgroup of $G_0$ constructed
in Case $1$, we may assume without loss of generality that $a$ takes values in
$S_0$. Let $S = {_aS_0}$ and consider the diagram
\[ \begin{array}{ccc}
H^1(K,S_0) & \xrightarrow{\pi_0} & H^1(K, G_0) \\
\uparrow \; f_S &      & \uparrow \; f_G       \\
H^1(K, S) & \xrightarrow{\pi} & \; H^1(K, G) \, . \end{array} \]
Here $f_S$ and $f_G$ are natural bijections. Since $\pi_0$ is surjective,
so is $\pi$.

\smallskip
{\it Case $3$.} Let $G$ be a connected reductive $k$-group. It is an
almost direct product of the semisimple $k$-group $H=[G,G]$ and the
central $k$-torus $C$ of $G$. Let $Z$ be the center of $H$.
Clearly, we have $C\cap H\le Z$. 
Consider the group $G' =G/Z$ and a natural morphism $f:G\to G'$.
By our construction, $G'$ is the direct product
of the torus $C/C\cap H$ and the adjoint group $H'=H/Z$.

Let $S'$ be the subgroup constructed in Case $2$ for $H'$ and let $\mu=
{_n(C/C\cap H)}$ be the $n$-torsion subgroup of the torus $C/C\cap H$,
where $n$ is the index of the minimal extension of $k$ splitting $C$.
Then for any extension $K/k$ a natural morphism $H^1(K,\mu\times S')
\to H^1(K,G')$ is surjective. We claim that $S=f^{-1}(\mu\times S')$
is as required, i.e. $H^1(K,S)\to H^1(K,G)$ is surjective.

Indeed, the exact sequences 
$1\to Z\to G\to G' \to 1$ and $1\to Z \to S\to S' \to 1$
give rise to a commutative diagram
$$
\begin{array}{ccccccc}
H^1(K,Z) & \xrightarrow{} & H^1(K,G) & \xrightarrow{g_1} & H^1(K, G')
& \xrightarrow{g_2} & H^2(K,Z) \\
\uparrow  & &\uparrow \, \pi &    & \uparrow \, \pi'   &  & \uparrow \, {id} \\
H^1(K,Z) & \xrightarrow{} & H^1(K, S) & \xrightarrow{h_1} &  H^1(K,\mu\times
S') &
\xrightarrow{h_2} & H^2(K,Z)
\end{array}
$$
Here $g_2,h_2$ are connecting homomorphisms.  Let $[a]\in H^1(K,G)$
and $[b]=g_1([a])$. Since $\pi'$ is surjective, there is a class $[c]\in
H^1(\mu\times S')$ such that $\pi'([c])=[b]$. Since $h_2([c])=
g_2 \pi'([c])=0$, there is $[d]\in H^1(K,S)$ such that $h_1([d])=[c]$.
Thus two classes $[a]$ and $\pi ([d])$ have the same image in $H^1(K,G')$.
By a twisting argument, one gets a surjective
map $H^1(K, {_dZ})\to g^{-1}(g_1([a]))$.
Since $Z\subset S$ and hence $_dZ\subset {_dS}$,
we have $[a]\in {\rm Im}\, \pi$. This completes the proof 
of Theorem~\ref{thm4}.
\qed

\smallskip

\begin{remark} \label{rem.constant} Our argument shows that if
$G$ is split and $k$ contains certain roots of unity,
then the subgroup $S$ in parts (b) and (c) can be taken to be 
a constant group.

More precisely, in part (c), $k$ needs to have a primitive 
root of unity of degree $n = |W(G_{\rm ss})|^2 \cdot |Z(G_{\rm ss})|$,
where $W(G_{\rm ss})$ and $Z(G_{\rm ss})$ denote, respectively,
the Weyl group and the center of the semisimple part $G_{\rm ss}$ of
$G$.

The same is true in part (b), except that $G$ needs to be replaced
by $G_{red} = G/R_u(G)$ in the above definition of $n$. 
\end{remark}

\section{Theorem~\ref{thm4} in the context of invariant theory}
\label{sect.geom}

For the rest of this paper $k$ will be 
an algebraically closed field of characteristic zero,
$K$ will be a finitely generated extension of $k$ and $G$ will
be a linear algebraic group defined over $k$.  
In this section we will introduce some terminology
in this context, discuss an invariant-theoretic 
interpretation of Theorem~\ref{thm4}(a) and use it to
give a simple proof of the no-name lemma. The third author 
would like to thank V. L. Popov for helpful suggestions
concerning this material.

\nsubsection{\bf $(G, S)$-sections.}
Recall that every
element of $H^1(K, G)$ is uniquely represented by a primitive
generically free $G$-{\em variety} $V$, up to birational
isomorphism.  That is, $k(V)^G = K$, the rational quotient
map $\pi \colon V \dasharrow V/G$ is a torsor over
the generic point of $V/G$, and this torsor is $\alpha$;
see~\cite[1.3]{popov}.  (Here ``$V$ is primitive" means 
that $G$ transitively permutes the irreducible 
components of $V$. In particular, if $G$
is connected then $V$ is irreducible.)

If $S$ is a closed subgroup of $G$ and $\alpha \in H^1(K, S)$
is represented by a generically free $S$-variety $V_0$, then
the image of $\alpha$ in $H^1(K, G)$ is represented by the
$G$-variety $G \ast_S V_0$, which is, by definition, the rational quotient
of $G \times V_0$ for the $S$-action given by
$s \colon (g, v_0) \mapsto (gs^{-1}, s\cdot v_0)$. We shall denote the image
of $(g, v_0)$ in this quotient by $[g, v_0]$.
Note that a rational quotient is, a priori, only defined 
up to birational isomorphism; however, a regular model 
for $G \ast_S V_0$ can be chosen so that
the $G$-action on $G\times V_0$ (by translations on the first factor)
descends to a regular $G$-action on $G\ast_S V_0$, making
the rational quotient map $G\times V_0 \dasharrow
G\ast_S V_0$ $G$-equivariant (via $g' \cdot [g, v_0] \mapsto [g' g , v_0]$);
see~\cite[2.12]{reichstein2}.
If $S$ is a finite group and $V_0$ is a quasi-projective $S$-variety
(which will be the case in the sequel) then we may take
$G \ast_S V_0$ to be the geometric quotient for the $S$-action on
$G \times V_0$, as in~\cite[Section 4.8]{pv}.

Now let $V$ be a $G$-variety. An $S$-invariant subvariety $V_0 \subset V$
is called a $(G, S)$-section if

\smallskip
(a) $G \cdot V_0$ is dense in $V$ and

\smallskip
(b) $V_0$ has a dense open $S$-invariant subvariety $U$ such that
$g \cdot u \in V_0$ for some $u \in U$ implies $g \in S$.

\smallskip
The above definition is due to Katsylo~\cite{katsylo}; 
sometimes a $(G, S)$-section is also called a {\em standard 
relative section} (see~\cite[1.7.6]{popov}) or {\em a relative 
section} with normalizer $S$ (see~\cite[Section 2.8]{pv}).
A $G$-variety $V$ is birationally isomorphic to
$G \ast_S V_0$ for some $S$-variety $V_0$ if and only if $V$ has
a $(G, S)$-section; see~\cite[Section 2.8]{pv}.
In this context Theorem~\ref{thm4}(a) can be rephrased as follows:

\smallskip
\noindent
{\bf Theorem~\ref{thm4}$'$}: {\em Every generically
free $G$-variety has a $(G, S)$-section, where $S$ is
a finite subgroup of $G$.}

\smallskip
Recall that a subvariety $V_0$ of a generically free $G$-variety
$V$ is called a {\em Galois quasisection} if the rational quotient map
$\pi \colon V \dasharrow V/G$  restricts to a dominant map
$V_0 \dasharrow V/G$, and the induced field extension $k(V_0)/k(V)^G$ is
Galois. If $V_0$ is a Galois quasisection then the finite group
$\Gamma(V_0) := \Gal(k(V_0)/k(V)^G)$ is called the Galois 
group of $V_0$; see~\cite{galitskii} or \cite[(1.1.1)]{popov}. 
(Note $\Gamma(V_0)$ is not required to be related to $G$ 
in any way.) The following theorem is due to 
Galitskii~\cite{galitskii}; cf. also~\cite[(1.6.2) and (1.17.6)]{popov}.

\begin{thm} \label{thm.galitskii} If $G$ is connected then every 
generically free $G$-variety has a Galois quasisection.
\end{thm}

A $(G, S)$-section is clearly a Galois quasisection with 
Galois group $S$. Hence, Theorem~\ref{thm4}$'$ 
(or equivalently, Theorem~\ref{thm4}(a)) may be viewed as
an extension of Theorem~\ref{thm.galitskii}.
Note that the Galois group $\Gamma(V_0)$ of the
Galois quasisection $V_0$ constructed in the proof 
of Theorem~\ref{thm.galitskii} is isomorphic 
to a subgroup of the Weyl group $W(G)$; cf. \cite[Remark 1.6.3]{popov}.
On the other hand, the group $S$ in our proof of
Theorem~\ref{thm4}(a), is an extension of $W(G)$ 
by a finite abelian group.  Enlarging the finite group $S$ 
may thus be viewed as ``the price to be paid" for 
a section with better properties.

\nsubsection{\bf The no-name lemma.}
A $G$-bundle $\pi \colon V \lra X$ is an algebraic vector bundle
with a $G$-action on $V$ and $X$ such that $\pi$ is $G$-equivariant
and $g$ restricts to a linear 
map $\pi^{-1}(x) \lra \pi^{-1}(gx)$ for every $x \in X$.

\begin{lem}[No-name Lemma] \label{no-name}
Let $\pi \colon V \lra X$ be a $G$-bundle of rank
$r$. Assume that the $G$-action on $X$ is generically free.
Then there exists a birational isomorphism 
$\pi \colon V \stackrel{\simeq}{\dasharrow}
X \times \bbA^r$ of $G$-varieties such that the following diagram commutes
\begin{equation} \label{e.no-name1}
\xymatrix{ V \ar@{-->}[r]^{\phi \quad } 
\ar@{->}[d]^{\pi} & X \times \bbA^r 
\ar@{->}[dl]^{\operatorname{pr}_1}  \cr X & } 
\end{equation}
Here $G$ is assumed to act trivially on $\bbA^r$, and
$\operatorname{pr}_1$ denotes the projection to the first factor.
In particular, $k(V)^G$ is rational over $k(X)^G$.
\end{lem}

The term ``no-name lemma", due to Dolgachev~\cite{dolgachev}, 
reflects the fact that
this result was independently discovered by many researchers.
In the case where $G$ is a finite group, Lemma~\ref{no-name} 
(otherwise known as Speiser's lemma, see~\cite{speiser})
may be viewed as a restatement of Hilbert's Theorem 90. 
In this case there are many proofs in the literature; 
see, e.g.,~\cite[Proposition 1.1]{em}, \cite[Proposition 1.3]{lenstra},
\cite[Appendix 3]{shafarevich} or \cite[Section 4]{cs}. 
For algebraic groups $G$ Lemma~\ref{no-name} was 
noticed more recently (the earliest reference we know is~\cite{bk}). 
This fact is now widely known and much used; however, as 
Colliot-Th\'el\`ene and Sansuc observed in~\cite[Section 4]{cs},
a detailed proof has never been published.  We will now use 
Theorem~\ref{thm4}(a) (or equivalently, Theorem~\ref{thm4}$'$
above) to give a simple argument reducing the general 
case of the no-name lemma to the case of a finite group.

\smallskip
{\em Proof of the no-name lemma:} By Theorem~\ref{thm4}$'$ 
$X$ has a $(G, S)$-section
$X_0$ for some finite subgroup $S$ of $G$. Then $V_0 = \pi^{-1}(X_0)$ is
a $(G, S)$-section for $V$; cf.~\cite[(1.7.7), Corollary 2]{popov}.
In other words, $X \simeq X_0 *_S G$ and
$V \simeq V_0 *_S G$, where $\simeq$ denotes birational 
isomorphism of $G$-varieties.

Note that $V_0$ is an $S$-vector bundle over $X_0$.
Since we know that the no-name lemma holds 
for $S$, there is an $S$-equivariant birational isomorphism 
$\phi_0 \colon V_0 \stackrel{\simeq}{\dasharrow} X_0 \times \bbA^r$ such that 
the diagram of $S$-varieties
\begin{equation} \label{e.no-name2} \xymatrix{
V_0 \ar@{-->}[r]^{\phi_0 \quad}  \ar@{->}[d]^{\pi} & X_0 \times \bbA^r 
\ar@{->}[dl]^{\operatorname{pr}_1}  \cr X_0 & } 
\end{equation}
commutes.  Taking the homogeneous fiber product of this diagram
with $G$, we obtain
\[ \xymatrix{
V \simeq V_0*_S G \ar@{-->}[r]^{\phi \quad \quad \quad} \ar@{->}[d]^{\pi} & 
(X_0 \times \bbA^r)*_S G  \simeq X \times \bbA^r
\ar@{->}[dl]^{\operatorname{pr}_1}  \cr X \simeq X_0 *_S G \, , & } 
\]
where $\phi = \phi_0 *_S G$.
\qed

\begin{remark}  The above argument can be naturally rephrased
in cohomological terms. Let $K = k(X)^G = k(X_0)^S$. Then  
Lemma~\ref{no-name} is equivalent to the following assertions:
(i) $k(V)^G = K(t_1, \dots, t_r)$
and (ii) $\alpha_V$ is the image of $\alpha_X$ under the 
restriction map $H^1(K, G) \lra H^1(K(t_1, \dots, t_r), G)$. 

Diagram~\eqref{e.no-name2} tells us that 
(i)$'$ $k(V_0)^S = K(t_1, \dots, t_r)$ and (ii)$'$ $\alpha_{V_0}$ is 
the image of $\alpha_{X_0}$ under the natural 
map $H^1(K, S) \lra H^1(K(t_1, \dots, t_r), S)$. (i) follows immediately 
from (i)$'$, and (ii) follows from (ii)$'$ by considering the natural diagram
\[ \xymatrix{\alpha_{X_0} \in H^1(K, S) \ar@{->}[r] \ar@{->}[d]^{\text{res}} &  
H^1(K(t_1, \dots, t_r), S) \ni \alpha_{V_0} 
 \ar@{->}[d]^{\text{res}}  \cr 
\alpha_X \in H^1(K, G) \ar@{->}[r] &  
H^1(K(t_1, \dots, t_r), G) \ni \alpha_V \, .} \]
\end{remark}

\section{Preliminaries on $G$-covers}
\label{sect.thm2a}

Let $G$ be a finite group.  We shall call 
a finite morphism $\pi \colon X' \lra X$ of algebraic 
varieties a $G$-{\em cover}, if $X$ is irreducible, 
$G$ acts on $X'$, so that $\pi$ maps every $G$-orbit in $X'$ to 
a single point in $X$, and $\pi$ is a $G$-torsor over 
a dense open subset $U$ of $X$. We will express 
the last condition by saying that $\pi$ is {\em unramified} 
over $U$. Restricting $\pi$ to the generic point of $X$, we obtain 
a torsor $\alpha \in H^1(k(X), G)$ over $\Spec \, k(X)$. 
In this situation we shall say that $\pi$ represents $\alpha$.
If a cover $\pi \colon X' \lra X$ is unramified over all of $X$, 
then we will simply say that $\pi$ is {\em unramified}. 

Recall that $\alpha \in H^1(K, G)$ is called  
{\em unramified} if it lies in the image of $H^1(R, G) \lra H^1(K, G)$
for every discrete valuation ring $k \subset R \subset K$ and
{\em projective}, if it is represented by 
an unramified $G$-cover $\pi \colon X' \lra X$
over a complete (or equivalently, projective) variety $X$.

\begin{lem} \label{lem.finite}
Let $G$ be a finite group, $K$ be a finitely
generated extension of an algebraically closed base field $k$
of characteristic zero, and $\alpha \in H^1(K, G)$.
Then the following assertions are equivalent:

\smallskip
(a) $\alpha$ is represented by a projective $G$-variety $V$
(in the sense of Section~\ref{sect.geom}), such that every 
element $1 \ne g \in G$ acts on $V$ without fixed points,

\smallskip
(b) $\alpha$ is projective, and

\smallskip
(c) $\alpha$ is unramified.
\end{lem}

\smallskip
Note that condition (b) can be rephrased by saying that $\alpha$
has trivial fixed point obstruction; see Section~\ref{sect.fpo}.

\smallskip
\begin{proof}
$(a) \Rightarrow (b)$: The $G$-action on $V$ has 
a geometric quotient $\pi \colon V \lra X$, where
$X$ is a projective variety; cf., e.g.,~\cite[Section 4.6]{pv}.
We claim that $\pi$ is a torsor over $X$. Indeed,
we can cover $V$ by $G$-invariant affine open 
subsets $V_i$. The quotient variety $X$ is then covered 
by affine open subsets $X_i = \pi(V_i)$, moreover, 
$\pi_i = \pi_{|V_i} \colon V_i \lra X_i$ is the geometric
quotient for the $G$-action on $V_i$; 
see~\cite[Theorem 4.16]{pv}. It is thus enough to show that
$\pi_i \colon V_i \lra X_i$ is a torsor for each $i$. This 
is an immediate corollary of the Luna Slice
Theorem; see, e.g.,~\cite[Theorem 6.1]{pv}.
       
\smallskip
(b) $\Rightarrow$ (c): Suppose $\alpha$ is 
represented by a $G$-torsor $V \lra X$, where $X$
is a projective variety with $k(X) = K$.
We want to prove that for any discrete valuation ring
$R \subset K$ the class $\alpha$ belongs to the image
$H^1(R,G)\to H^1(K,G)$.

Indeed, the ring $R$ dominates a point in $X$;
denote this point by $D$.  
Consider the canonical map $\Spec \, R \to X$
sending the closed point in $\Spec \, R$ to $D$ and the
generic point of $\Spec \, R$ into the generic point of $X$.
Take the fiber product $(\Spec \,  R) \times _X V$. It follows
immediately from this construction that the $G$-torsor
$$
(\Spec \, R) \times _X V \to \Spec \, R
$$
is as required, i.e. its image under
the map $H^1(R,G)\to H^1(K,G)$ is $\alpha$.

\smallskip
(c) $\Rightarrow$ (a): Let $V$ be a smooth projective $G$-variety
representing $\alpha$ and let $\pi \colon V \lra X$ be 
the geometric quotient.  Note that $X$ is normal.
We want to show that every $1 \ne g \in G$
acts on $V$ without fixed points. Assume the contrary:  
$g v = v$ for some $v \in V$.
By~\cite[Theorem 9.3]{ry2} (with $s = 1$ and
$H_1 = \mathopen< g \mathclose>$), after performing
a sequence of blowups with smooth $G$-invariant centers on $V$,
we may assume that the fixed point locus $V^g$ of $g$
contains a divisor $D \subset V$. If $R = \mathcal{O}_{X, \pi(D)}$ 
is the local ring of the divisor $\pi(D)$ in $X$ then
$\alpha$ does not lie in the image of the natural morphism 
$H^1(R, G) \lra H^1(K, G)$, a contradiction.
\end{proof}

\begin{remark} \label{rem.finite} Our proof of the implication
(b) $\Rightarrow$ (c) does not use the fact that $G$ is a finite group. 
This implication is valid for every linear algebraic group $G$.
\end{remark} 

\section{Proof of Theorem~\ref{thm2}}
\label{sect.thm2b}

Let $S$ be the finite subgroup of $G$ given by Theorem~\ref{thm4}(a).
Then $\alpha \in H^1(K, G)$ is the image of some $\beta \in H^1(K, S)$. 
Examining the diagram
\[ \xymatrix{H^1(X, S) \ar@{->}[r] \ar@{->}[d] &  H^1(L, S) \ni \beta_L 
 \ar@{->}[d]  \cr 
H^1(X, G) \ar@{->}[r] &  H^1(L, G) \ni \alpha_L \, ,} \] 
%
where $X$ is a complete variety and $L = k(X)$,
we see that if Theorem~\ref{thm2} holds for $S$ then it holds for $G$.

 From now on we may assume that $G$ is a finite group.
In this case Theorem~\ref{thm2} can be restated as follows. 

\begin{prop} \label{prop.aby2}
Let $G$ be a finite group, $k$ be an algebraically closed
base field of characteristic zero, $K/k$ be a finitely 
generated extension, and $\alpha \in H^1(K, G)$. Then there exists an  
abelian field extension $L/K$ such that $\alpha_L$ is represented by
an unramified $G$-cover $\pi \colon Z' \lra Z$, where $Z$ and $Z'$
are projective varieties.
\end{prop}

The rest of this section will be devoted to proving
Proposition~\ref{prop.aby2}.  We begin with the following lemma.

\begin{lem} \label{lem.aby1}
Let $G$ be a finite group. Then every $\alpha \in H^1(K, G)$
is represented by a $G$-cover $\pi \colon X' \lra X$ such that

\smallskip
(a) $X'$ is normal and projective,

\smallskip
(b) $X$ is smooth and projective,

\smallskip
(c) there exists a normal crossing divisor $D$ on $X$ such that
$\pi$ is unramified over $X - D$.
\end{lem}

\begin{proof} Suppose $\alpha$ is represented by a $G$-Galois
algebra $K'/K$. We may assume without loss of generality that
$K'$ is a field. Indeed, otherwise $\alpha$ is the image of some
$\alpha_0 \in H^1(K, G_0)$, where $G_0$ is a proper subgroup of $G$,
and we can replace $G$ by $G_0$ and $\alpha$ by $\alpha_0$.

Choose a smooth projective model $Y/k$ for $K/k$ and
let $\phi \colon Y' \lra Y$ be the normalization
of $Y$ in $K'$. Then $Y'$ is projective
(see~\cite[Theorem III.8.4, p. 280]{mumford}),
and by uniqueness of normalization
(see~\cite[Theorem III.8.3, pp. 277 - 278]{mumford}),
$G$ acts on $Y'$ by regular morphisms, so that $k(Y')$ is
isomorphic to $K'$ as a $G$-field (see~\cite[pp. 277 - 278]{mumford}).
We have thus shown that $\alpha$ can be represented
by a cover $\phi \colon Y' \lra Y$ 
satisfying conditions (a) and (b). We will now birationally modify 
this cover to obtain another cover $\pi \colon X' \lra X$ 
which satisfies condition (c) as well.

The cover $\phi$ is unramified over a dense open subset
of $Y$; denote this subset by $U$. Set
$E = Y - U$, and resolve $E$ to a normal crossing
divisor $D$ via a birational morphism  $\gamma \colon X \lra Y$.
Now consider the diagram
\[ \begin{array}{ccc}
X' & \lra & Y'  \\
\quad \downarrow \; \pi &   & \quad \downarrow \; \phi  \\
X & \xrightarrow{\; \; \gamma \; \;} & \; Y \, , \end{array} \]
where $X'$ is the normalization of $X$ in $K'$.
By our construction $X$ is smooth and $X'$ is normal. Moreover,
since $\gamma$ is an isomorphism over $U$,
$\pi$ is unramified over $X - D = \phi^{-1}(U)$,
as desired.
\end{proof}

We are now ready to complete the proof of
of Proposition~\ref{prop.aby2}. Our argument will be based
on~\cite[Theorem 2.3.2]{gm}, otherwise known
as ``Abhyankar's Lemma", which describes the local
structure of a covering, satisfying conditions
(a) - (c) of Lemma~\ref{lem.aby1}, in the \'{e}tale topology.
We thank K. Karu for bringing this result to our 
attention.

\smallskip
Let $\pi \colon X' \lra X$ be a $G$-cover of projective varieties
representing $\alpha$
and satisfying conditions (a) - (c) of Lemma~\ref{lem.aby1}.
Denote the irreducible components of $D$ by $D_1, \dots, D_s$.

Since $X$ is smooth, each $x \in X$ has an affine open neighborhood
$U_x$ where each $D_j$ is principal, i.e., is given by
$\{ a_{x, j} = 0 \}$ for some $a_{x, j} \in {\cal O}_X (U_x)$ (possibly
$a_{x, j} = 1$ for some $x$ and $j$).
By quasi-compactness, finitely many of these open
subsets, say, $U_{x_1}, \dots, U_{x_n}$ cover $X$. To simplify our
notation, we set $U_i = U_{x_i}$ and $a_{i j} = a_{x_i, j}$.

Now let $b_{ij}$ be an $|G|$th root of $a_{ij}$ in
the algebraic closure of $K = k(X)$ and $L = K(b_{ij})$,
where $i$ ranges from $1$ to $n$ and $j$ ranges from $1$ to $s$.
Suppose $\gamma \colon Z \lra X$ is the normalization of $X$
in $L$ and $Z' = X' \times_{X} Z$.  Since we are
assuming that $k$ is algebraically closed of characteristic zero
(and in particular, $k$ contains a primitive $|G|$th
root of unity), $L/K$ is an abelian extension.
It is also easy to see from our construction
that $Z$ and $Z'$ are projective, $Z$ is normal, and
the natural projection $\pi' \colon Z' \lra Z$ is a $G$-cover,
which represents $\alpha_L \in H^1(L, G)$. To sum up,
we have constructed the following diagram of morphisms:
\[ \begin{array}{ccc}
          Z'  &  \rightarrow & X' \\
          \, \quad \downarrow \; \psi &      & \quad \downarrow \; \pi \\
          Z &  \xrightarrow{\gamma} & \; X \, .
\end{array} \]
It remains to show that the $G$-cover $\psi$ is unramified.
Suppose we want to show that $\psi$ is unramified at $z_0 \in Z$.
Since the open sets $U_1, \dots, U_n$ cover $X$,
$x_0 = \gamma (z_0)$ lies in $U_i$ for some $i = 1, \dots, n$.
By Abhyankar's lemma~\cite[Theorem 2.3.2]{gm}, there
exists an abelian subgroup $H \simeq \bbZ/n_1 \times \dots
\times \bbZ/n_s \bbZ$ of $G$ (possibly with $n_j = 1$ for some $j$)
and a (Kummer) $H$-Galois cover
\[ V_i = \{ (x, t_1, \dots, t_s) \; | \; t_1^{n_1} = a_{i, 1} \, , \; \dots,
\; t_s^{n_s} = a_{i, s} \} \subset  U_i \times \bbA^s \, , \]
such that the $G$-covers $\pi \colon X' \lra X$ and
$\phi \colon G *_H V_i \lra U_i$ are isomorphic over an \'{e}tale
neighborhood of $x_0$ in $X$.  (Here the natural projection $V_i \lra U_i$
is an $H$-cover, and $G *_H V_i \lra U_i$ is the $G$-cover induced 
from it; for a definition of $G *_H V_i$, see Section~\ref{sect.geom}.)

Now recall that by our construction the elements 
$b_{ij} \in L = k(Z)$ satisfy $b_{ij}^{|G|} = a_{ij} 
\in {\cal O}_X(U_i)$. In particular, they are integral 
over $U_i$ and thus they are regular function on $\gamma^{-1}(U_i)$. 
Since $n_j$ divides $|G|$ for every $j = 1, \dots, s$,
the pull-back of $\phi$ to $Z$ splits over an \'{e}tale neighborhood
of $z_0$; hence, so does $\psi$ = pull-back of $\pi$.
In other words, $\psi$ is unramified at $z_0$, as claimed.
This completes the proof of Proposition~\ref{prop.aby2}.
\qed

\section{An example}
\label{sect.g2}

It is well known that there exist non-trivial projective elements
in $H^1(K, \PGLn)$ for every $n \ge 2$ (for suitable $K$).
In this section we use a variant of a construction 
of Colliot-Th\'el\`ene and Gabber~\cite{cg} to show that,
for certain $K$, such elements exist in $H^1(K, G_2)$ as well.

\begin{prop} \label{prop.g2}
Let $k$ be an algebraically closed base field 
of characteristic zero such that $\trdeg_{\bbQ}(k) \ge 3$.
(Note that the last condition is satisfied by every 
uncountable field.) Then there exist a smooth 
projective 3-fold $X/k$ with function field $K = k(X)$
and a projective non-trivial class $\alpha \in H^1(K, G_2)$.
\end{prop}

Note that no such examples can exist if $X$ is a curve
or a surface, since in this case $H^1(k(X), G_2) = \{ 1 \}$;
see~\cite{bp}.  

\smallskip
\begin{proof} 
Let $E_1$, $E_2$, $E_3$ be elliptic curves. For $i = 1, 2, 3$
choose $p_i, q_i \in E_i$
so that $p_i \ominus q_i$ is a point of order $2$. (Here $\ominus$ denotes 
denotes subtraction with respect to the group operation on $E_i$.)
Then $2p_i - 2q_i$ is a principal divisor on $E_i$ and $p_i - q_i$ is not; 
see, e.g., \cite[Corollary 3.5]{silverman}. 
Thus $2p_i - 2q_i = \operatorname{div}(f_i)$, where $f_i \ne 0$ is 
a rational function on $E_i$, which is not a complete square.
Adjoining $\sqrt{f_i}$ to $k(E_i)$, we obtain an irreducible unramified
$\bbZ/ 2\bbZ$-cover $\pi_i \colon E_i' \lra E_i$. (Note that
by the Hurwitz formula, $E_i'$ is also an elliptic curve.)

Now set $X = E_1 \times E_2 \times E_3$ and $K = k(X)$,
$S = (\bbZ/2 \bbZ)^3$, and consider the element 
$\beta \in H^1(k(X), S)$, represented 
by the $S$-cover
\[ \pi = (\pi_1, \pi_2, \pi_3) \colon 
E_1' \times E_2' \times E_3' \lra E_1 \times E_2 \times E_3 = X \, . \]
Since $\pi$ is an unramified cover, $\beta$ is projective.

We now recall that the exceptional group $G_2/k$ contains a
unique (up to conjugacy) maximal elementary abelian $2$-group 
$i \colon S= (\bbZ / 2 \bbZ)^3 \hookrightarrow G_2$. Set
$\alpha = i_*(\beta) \in H^1(K, G_2)$. Since $\beta$ is 
projective, so is $\alpha$.  It thus remains to show that 
$\alpha \ne 1$ in  $H^1(K, G_2)$ (for a suitable choice 
of $E_i$ and $E_i'$).

The cohomology set $H^1(K,G_2)$ classifies octonion algebras 
or equivalently, $3$-fold Pfister forms; cf.~\cite[Theorem 9]{serrepp}.
By~\cite[\S 22.10]{gms}, the map
$$ H^1(K, S) = \Bigl(K^\times / (K^\times)^2 \Bigr)^3  
\stackrel{i_*}{\lra} H^1(K, G_2) $$
is non-trivial; hence, it sends $(a_1, a_2, a_3) \in
\Bigl(K^\times / (K^\times)^2 \Bigr)^3$ to 
the class of the $3$-Pfister form
$\langle \langle a_1, a_2, a_3 \rangle\rangle$; see 
\cite[Theorem 27.15]{gms}. By our construction, 
$\beta \in H^1(K, S)$ corresponds to $(f_1, f_2, f_3) \in 
\Bigl(K^\times / (K^\times)^2 \Bigr)^3$. Thus $\alpha = i_*(\beta)$
is non-split in $H^1(K, G_2)$ if and only if the 3-fold Pfister form
$\langle \langle f_1, f_2, f_3 \rangle\rangle$ is nonsplit or,
equivalently, if $(f_1) \cup (f_2) \cup (f_3)  \neq 0$ in 
$H^3(k(X),\bbZ / 2 \bbZ)$; see~\cite[Corollary 3.3]{elman-lam}. 

Since we are assuming that $\trdeg_{\bbQ}(k) \ge 3$, we can choose
elliptic curves $E_1$, $E_2$ and $E_3$ so that their
$j$-invariants  are algebraically independent over $\mathbb{Q}$.
We now appeal to a theorem of Gabber
(\cite[p. 144]{cg}), which says that  
$(f_1) \cup (f_2) \cup (f_3)  \neq 0$ in $H^3(k(X), \bbZ / 2 \bbZ)$. 
Hence, $\alpha \ne 1$ in $H^1(K, G)$, as claimed. This completes the proof
of Proposition~\ref{prop.g2}.
\end{proof}

%
%

\section{The fixed point obstruction}
\label{sect.fpo}

We now recall the notion of {\em fixed point obstruction}
from~\cite[Introduction]{ry8}.
Suppose $\alpha \in H^1(K, G)$ is represented by 
a generically free primitive $G$-variety $V$ (as in
Section~\ref{sect.geom}).
We shall say that a subgroup of $G$ is {\em toral} if it
lies in a subtorus of $G$ and {\em non-toral} otherwise.
If $V$ (or any $G$-variety birationally isomorphic to it) has a smooth
point fixed by a non-toral diagonalizable subgroup $H \subset G$,
then we shall say that $V$ (or equivalently, $\alpha$) has
{\em non-trivial fixed point obstruction};  
cf.~\cite[Introduction]{ry8}.  Note that after birationally
modifying $V$, we may assume that $V$ is smooth and complete
(or even projective, see, e.g.,~\cite[Proposition 2.2]{ry2}),
and that the fixed point obstruction can be detected 
on any such model. In other words, if $V$ and $V'$ are 
smooth complete birationally isomorphic $G$-varieties
then $V^H = \emptyset$ if and only if $(V')^H = \emptyset$ for
any diagonalizable subgroup $H \subset G$;
see~\cite[Proposition A2]{ry}. If $V^H = \emptyset$ 
for every diagonalizable non-toral subgroup $H \subset G$ 
(and $V$ is smooth and complete),
then we will say that $V$, or equivalently $\alpha$, 
has {\em trivial fixed point obstruction}.  

If $\alpha$ is split (i.e., $\alpha = 1$ in $H^1(K, G)$)
then by~\cite[Lemma 4.3]{ry2} $\alpha$ has trivial fixed 
point obstruction. We will now extend this result as follows.

\begin{prop} \label{prop.fpo}
If $\alpha \in H^1(K, G)$ is projective
then $\alpha$ has trivial fixed point obstruction.
\end{prop}

\begin{proof} Let $\overline{G}$ be a smooth projective 
$G \times G$-variety, which contains $G$ as a dense open orbit.
(Here we are viewing $G$ as a $G \times G$-variety
with respect to left and right multiplication). To construct $\oG$,
we use a theorem of Kambayashi, which says that $G$ can be
$G \times G$-equivariantly embedded into $\bbP(V)$ for some
linear representation $G \times G \lra \GL(V)$; see~\cite[Theorem 1.7]{pv}.
Taking the closure of $G$ in $\bbP(V)$, and 
$G \times G$-equivariantly resolving its singularities, 
we obtain $\overline{G}$ with desired properties.  

For $\overline{g} \in \overline{G}$, we will write 
$g_1 \cdot \overline{g} \cdot g_2^{-1}$
instead of $(g_1, g_2) \cdot \overline{g}$; the reason for this notation
is that for $\overline{g} \in G$, 
$(g_1, g_2) \cdot \overline{g} = g_1 \overline{g} g_2^{-1} \in G$.

Since $\alpha$ is projective, it can be represented
by a $G$-torsor $\pi \colon Z \lra X$ over a smooth projective 
irreducible variety $X$. (Here $K = k(X)$.) 
We will now construct a smooth complete $G$-variety 
$\overline{Z}$ representing $\alpha$ (i.e., birationally 
isomorphic to $Z$) by "enlarging" each fiber of $\pi$ from $G$
to $\oG$.

Let $U_i \to X$, $i \in I$ be an \'{e}tale covering 
which trivializes $\pi$. Then $\pi$ is described by 
the transition maps 
$f_{ij} \colon U_{ij} \times G \lra U_{ij} \times G$ on 
the pairwise ``overlaps" $U_{ij}$; here each $f_{ij}$ is an automorphism
of the trivial $G$-torsor $U_{ij} \times G$ on $U_{ij}$. 
($G$ acts trivially on $U_{ij}$ and by left translations 
on itself.) These transition maps satisfy 
a cocycle condition (for Cech cohomology)
which expresses the fact that they are compatible 
on triple ``overlaps" $U_{hij}$. It is easy to see that 
$f_{ij}$ is given by the formula
\begin{equation} \label{e.h}
f_{ij} (u, g) = (u, g \cdot h_{ij}(u)) \, ,  
\end{equation}
for some morphism $h_{ij} \colon U_{ij} \lra G$.
(In fact, $h_{ij}(u) = \pr_2 \circ f_{ij}(u, 1_G)$, where $\pr_2 \colon
U_{ij} \times G \lra G$ is the projection to the second factor.) 
Formula~\eqref{e.h} can now be used to extend 
$f_{ij}$ to a $G$-equivariant automorphism
\[ \overline{f_{ij}} \colon U_{ij} \times \oG \lra U_{ij} \times \oG \, , \]
where $G$ acts on $\oG$ on the left. Since $f_{ij}$ satisfies 
the cocycle condition and $G$ is dense in $\overline{G}$, we 
conclude that $\overline{f_{ij}}$ satisfy the cocycle condition 
as well. By descent theory, the transition maps $\overline{f_{ij}}$
patch together to yield a variety $\overline{Z}$ and a commutative
diagram of morphisms
\[   \begin{array}{ccccc}  Z &   & \hookrightarrow &  & \overline{Z} \\
                 & \pi \;  \searrow &   & \swarrow \; \overline{\pi}   &        \\
                 &     &   X    &    &  
\end{array} \]
which locally (in the \'{e}tale topology) looks like
\[   \begin{array}{ccccc}  U_i \times G &   & \hookrightarrow & 
                                                   & U_i \times \oG \\
                 & \pi \;  \searrow &   & \swarrow \; \overline{\pi}   &        \\
                 &     &   U_i    &    &  
\end{array} \]
(The maps $\pi$ and $\overline{\pi}$ in the second diagram are
projections to the first component.) 
It is now easy to see that $\overline{Z}$ is smooth and proper 
over $X$ and $Z \hookrightarrow \overline{Z}$ is a $G$-equivariant 
open embedding. Indeed, these properties can be checked locally
(in the \'{e}tale topology) on $X$, where they are immediate from the
second diagram. Note also that since $\overline{Z}$ is proper
over $X$, and $X$ is projective over $k$, $\overline{Z}$ is complete 
as a $k$-variety. 

Having constructed a smooth complete model $\overline{Z}$ for 
$\alpha$, we are now ready to show that $\alpha$ has trivial 
fixed point obstruction. Suppose a diagonalizable
subgroup $H$ of $G$ has a fixed point in $z \in \overline{Z}$. We want 
to show that $H$ is toral in $G$.  Indeed, let
$F$ be the fiber of $\overline{\pi}$ containing $z$. 
By our construction $F \simeq \oG$ as $G$-varieties 
(here $\oG$ is viewed as a $G$-variety with respect to the left 
$G$-action). We conclude that $H$ has a fixed point in $\oG$. 
Since $\oG$ has $G$ as a $G$-invariant dense open subset, 
it is split as a $G$-variety
(i.e., it represents the trivial class in $H^1(k, G)$),
\cite[Lemma 4.3]{ry2} now tells us that
$H$ is toral. This shows that
$\alpha$ has trivial fixed point obstruction, thus completing
the proof of Proposition~\ref{prop.fpo}.
\end{proof}

\begin{remark} \label{rem.fpo1} 
The fact that $G$ acts on $\oG$ both on the right and 
on the left was crucial in the construction of $\overline{Z}$ in
the above proof.  The action on the right was used to glue the transition
maps $\overline{f_{i, j}}$ together, and the action on the left
to define a $G$-action on $\overline{Z}$. If $G$ could only act on
$\oG$ on one side, we would still be able to construct $\overline{Z}$
as a variety; however, we would no longer be able to define a $G$-action
on it, extending the $G$-action on $Z$.
\end{remark}

\begin{cor} \label{cor.fpo2} 
There exist non-split elements $\alpha_n \in H^1(K_n, \PGL_n)$
($n = 2, 3, \dots$) and $\beta \in H^1(K, G_2)$ with trivial 
fixed point obstruction, for some finitely generated field
extensions $K_n/k$ and $K/k$.
\end{cor}

\begin{proof} Choose $\alpha_n$ and $\beta$ so that 
they are non-split and projective; cf. Section~\ref{sect.g2}.
\end{proof}

\begin{remark} \label{rem.fpo3}
By~\cite[Theorem 4]{ry8} for every prime number $p$
there exists a non-split $\alpha \in H^1(K, \PGL_p)$ 
such that $K$ is a purely transcendental extension of $k$ 
and $\alpha$ has trivial fixed point obstruction.
Such $\alpha$ are necessarily ramified 
and hence, cannot be projective. Thus 
the converse to Proposition~\ref{prop.fpo} is false.
\end{remark}

\section{Problem~\ref{prob.tits} and Hilbert's 13th problem}
\label{sect.hilbert}

\nsubsection{\bf An algebraic variant of Hilbert's 13th problem.}
Hilbert's 13th problem asks, loosely speaking, which continuous
functions in $n$ variables can be expressed as compositions 
of functions in $n-1$ variables. In this form the problem
was settled by Arnold~\cite{arnold} and Kolmogorov~\cite{kolmogorov},
who showed that any continuous function in $n$ variables can be expressed
as a composition of continuous functions in one variable and 
the addition function $f(x, y) = x + y$.
The algebraic variant of Hilbert's 13th problem, 
where ``continuous functions" are replaced by 
``algebraic functions", remains open. In modern 
language the problem can be stated as follows; cf.~\cite{as} 
or~\cite{dixmier}.

Let $E/F$ be a finite separable field extension (or, more generally, 
an \'{e}tale algebra) and assume that $F$ contains a copy 
of the base field $k$. Then the {\em essential dimension} 
$\ed_k(E/F)$ (or simply $\ed(E/F)$, if the reference to $k$ is 
clear from the context) is the minimal value of $\trdeg_k(F_0)$,
where the minimum is taken over all elements $a \in E$ and over all 
subfields $k \subset F_0 \subset F$ such that $E = F(a)$ and
$F_0$ contains every coefficient of the characteristic
polynomial of $a$; cf.~\cite{br1}, \cite{br2}.
For example, if $E/F$ is a non-trivial cyclic extension
of degree $n$ and $k$ contains a primitive $n$th root 
of unity then $\ed_k(E/F) = 1$, since in this case 
we can choose $a$ so that $a^n \in F$. Note also that
$\ed(E^{\#}/F) = \ed(E/F)$, where
$E^{\#}$ is the normal closure of $E$ over $F$; 
cf.~\cite[Lemma 2.3]{br1}.

We will now say that $E/F$ has level $\le d$ if
there exists a tower of finite field extensions
\begin{equation} \label{e.tower}
F = F_0 \subset F_1 \subset \dots \subset F_r
\end{equation}
such that $F \subset E \subset F_r$ and
$\ed_k(F_{i}/F_{i-1}) \le d$ for 
every $i = 1, \dots, r$. For example, if $k$ contains 
all roots of unity then every solvable extension 
$E/F$ has level $\le 1$ (because we can take~\eqref{e.tower}
to be a tower of cyclic extensions). The algebraic form 
of Hilbert's 13th problem then asks for the smallest integer
$s(n)$ such that the level of every degree $n$ extension $E/F$
is $\le s(n)$. (Here we are assuming that the base field $k$ 
is fixed throughout.) Not much is known about $s(n)$
(see~\cite{dixmier}); in particular, it is not known
if $s(n) > 1$ for any $n$. It is thus natural to ask, if, 
perhaps, $s(n) = 1$ for all $n$; this equality may be viewed
as an algebraic analogue of the above-mentioned 
theorem of Arnold and Kolmogorov; 
cf.~\cite[p. 90]{dixmier}. In fact, in the absence of evidence 
to the contrary, one can even ask for a particularly nice 
tower~\eqref{e.tower}, showing that $s(n) = 1$, namely for
a tower~\eqref{e.tower}, where $F_{r-1}/F$ is solvable 
(or even abelian) and $F_r/F_{r-1}$ has essential dimension 
$1$.  Equivalently,

 \begin{problem} \label{prob.hilbert}
 Let $k$ be an algebraically closed field of characteristic zero,
 $S$ be a finite group and $K/k$ be a field extension.
 Is it true that for every $\alpha \in H^1(K, S)$ there exists
 (i) an abelian extension $L/K$ such that $\ed(\alpha_L) \le 1$?
 or (ii) a solvable extension $L/K$ such that $\ed(\alpha_L) \le 1$?
 \end{problem}
 
Here $\alpha_L$ is represented by an $S$-Galois algebra $L'/L$ and
$\ed(\alpha_L)$ denotes the essential dimension of $L'/L$. Equivalently,
$\ed(\alpha_L)$ is the minimal value of $\trdeg_k(L_0)$ such that
$\alpha_L$ lies in the image of the natural map $H^1(L_0, S) \lra H^1(L, S)$
for some intermediate field $k \subset L_0 \subset L$.
(Note, that, since the base field is assumed to be algebraically closed,
$\ed(\alpha_L) = 0$ if and only if $\alpha_L$ is split.)
 
We do not know whether or not the assertions of
Problem~\ref{prob.hilbert} are true 
(cf. Remark~\ref{rem.abhyankar} below).  However, using
Theorem~\ref{thm4} we will show that, if true,
they have some remarkable consequences.
 
 \begin{thm} \label{main}
 Let $k$ be an algebraically closed field of characteristic zero.
 and let $K/k$ be a field extension. Denote the maximal abelian
 and the maximal solvable extensions of $K$ by
 $K_{ab}$ and $K_{sol}$ respectively.
 
 \smallskip
 (i) If Problem~\ref{prob.hilbert}(i) has an affirmative answer 
 then $\cd(K_{ab}) \le 1$.
 
 \smallskip
 (ii) If Problem~\ref{prob.hilbert}(ii) has an affirmative answer 
 then $\cd(K_{sol}) \le 1$.
 \end{thm}
 
\begin{remark} \label{rem.main}
The inequality $\cd(K_{ab}) \le 1$ is only known
in a few cases; in particular, for $K$ = a number field, 
or $K$ = a $p$-adic field by class field theory and
for $K = \mathbb{C}((X))((Y))$ by a theorem of
Colliot-Th\'el\`ene, Parimala and Ojanguren \cite[Theorem 2.2]{CTOP}.
If it were established, it would immediately imply an affirmative 
answer to Problem~\ref{prob.tits}. Another important consequence
would be a conjecture of Bogomolov~\cite[Conjecture 2]{Bog},
which asserts that $\cd \bigl(K^{(p)}_{ab} \bigr) \leq 1$,
where $K^{(p)}$ is a maximal prime-to-$p$ extension of $K$.
On the other hand, an affirmative
answer to Problem~\ref{prob.hilbert}(ii) would imply 
that ${\rm cd}_p\bigl( K(p) \bigr) \leq 1$, where $p$ 
is a prime number and $K(p)$ is the $p$-closure (i.e the maximal 
$p$--solvable extension) of $K$, thus giving an affirmative answer to
a question of J. K\"onigsmann; cf.~\cite[Question 5.3]{Koe}.
\end{remark}

        
\begin{remark} \label{rem.abhyankar} The third author would like
to take this opportunity to correct a misstatement he made in
\cite[Introduction]{br1}. The identity $d'(6) = 2$, which is
attributed to Abhyankar~\cite{abhyankar} at the bottom
of p. 161 in \cite{br1}, would, if true, give a negative
answer to Problem~\ref{prob.hilbert}(ii) for the
symmetric group $G = \Sym_6$. In fact, the version of Hilbert's 13th
problem considered in~\cite{abhyankar} is quite different from
ours; the base extensions that are allowed there are integral
ring extensions, rather than field extensions. For this reason
the identity $d'(6) = 2$ does not follow from the results
of~\cite{abhyankar} and, to the best of our knowledge,
Problem~\ref{prob.hilbert} is still open, even in the case where
$S$ is the symmetric group $\Sym_6$.
\end{remark}
        
\nsubsection{\bf Proof of Theorem~\ref{main}.}
We begin with some preliminary facts.
Recall that a field $F$ has cohomological dimension $\leq 1$
if and only if the Brauer group $\Br(F')$ is trivial for any
separable finite field
extension $F'/F$; see~\cite[Proposition II.3.5]{serre-gc}.
It will be convenient for us to work with \'etale $K$-algebras,
rather than just separable field extensions of $K$.
Recall that a $K$-\'etale algebra is
a finite product $E = K_1 \times K_2 \times \cdots \times
K_n$ of finite separable extensions $K_i/K$.
The Brauer group of $E$
is  $\Br(E) = \oplus_i \Br(K_i)$; an element of this group
is represented by an $n$-tuple ${\cal A}=({\cal A}_i/K_,i)_{i=1,..,n}$ of
central simple algebras. Note that ${\cal A}$ is an Azumaya
algebra over $E$.  
Given a field $F$, we have
\begin{equation} \label{e.brauer}
\cd(F) \le 1 \quad \Longleftrightarrow \quad
\text{$\Br(E)=0$ for any \'etale algebra $E/F$;}
\end{equation}
see~\cite[Proof of Theorem III.2.2.1]{serre-gc} or~\cite[Lemma 10.11]{fj}.  

\begin{lem} \label{lem.brauer} The following are equivalent:

\smallskip
(a) $cd(K_{ab} )\leq 1$,

\smallskip

(b) For any \'etale algebra $E/K$, the restriction map
$\Br(E) \lra {\Br(E \otimes_K K_{ab})}$ is trivial.

\smallskip
\noindent
Moreover, the lemma remains true if $K_{ab}$ is replaced by $K_{sol}$.
\end{lem}

\begin{proof}
(a) $\Rightarrow$ (b): immediate from~\eqref{e.brauer}.

\smallskip
(b) $\Rightarrow$ (a): Let $B/K_{ab}$ be an \'etale algebra.
There exists a finite abelian subextension $K'/K$ of
$K_{ab}/K$ and an \'etale
algebra $B'/K'$ such that $B' \otimes_{K'} K_{ab}=B$.
We have
$$
B = \limind\limits_{K' \subset L \subset K_{ab}} B' \otimes_{K'} L,
$$
where the limit is taken on subfields $L$ of $K_{ab}$ finite over $K'$.
Consequently,
$$
\Br(B) = \limind\limits_{K' \subset L \subset K_{ab}} \Br(B' \otimes_{K'} L),
$$
and (b) implies that $\Br(B)=0$. (a) now follows from~\eqref{e.brauer}.

The proof remains unchanged if $K_{ab}$ is replaced by $K_{sol}$.
\end{proof}

We are now ready to proceed with the proof of Theorem~\ref{main}(i).
We start with the group $G=(\PGLn)^m \sdp \Sym_m$.
By Theorem \ref{thm4}(a),  there exists a finite subgroup $S$
of $G$ such that the natural
homomorphism $H^1(K, S) \lra H^1(K,G)$ is surjective
The group $\Sym_m$ is the automorphism group
of the trivial \'etale algebra, so by Galois descent
the set $H^1(K, \Sym_m)$ classifies $m$-dimensional
\'etale algebras. By~\cite[Corollary I.5.4.2]{serre-gc},
the fiber of the map $H^1(K,G) \lra H^1(K, \Sym_m)$
at $[E] \in H^1(k, \Sym_m)$ is
$$
H^1\Bigl(K, _E(\PGLn^m) \Bigr) /  _E(\Sym_m) \, ,
$$
with $_E(\PGLn^m)$ and $_E(\Sym_m)$ are the twisted groups by
the \'etale algebra $E/K$. Since $G\to S_m$
has a section, the map $_EG(K)\to {_E}(S_m)(K)$ is surjective.
Then $_E(S_m)$ acts trivially on $H^1(K,{_E}(\PGLn^m))$
and hence the fiber at $[E]$ is  $H^1(K,{_E}(\PGLn^m))$.
By definition of the Weil restriction,
we have $_E (\PGL_n^m)=R_{E/k}(\PGL_n)$.
We identify
$H^1\Bigl(K, _E (\PGL_n^m) \Bigr) =
H^1(E, \PGL_n)$ by the Shapiro isomorphism. Thus
$$
H^1(K,G)= \bigsqcup\limits_{[E] \in H^1(K,\Sym_m)} H^1(E,PGL_n)\,.
$$
An element of $H^1(K,G)$ is then given
by an Azumaya algebra ${\cal A}/E$ of degree $n$
defined over a $K$-\'etale algebra $E$ of rank $m$.
By Theorem~\ref{thm4}(a), every
class $[{\cal A}/E]$ comes from a class $\alpha \in H^1(K,S)$.

We now apply the assertion of Problem~\ref{prob.hilbert}(i) 
to the group $S$ and the class $\alpha$. There exists an
abelian extension $L/K$, a $k$-curve $C$ and a map
$k(C) \subset L$ such that  the restriction of the
class $\alpha$ in $H^1(L,S)$  belongs
to the image of $H^1(k(C),S) \lra H^1(L,S)$.
The commutative diagram of restriction maps
\[ \begin{array}{ccccc}
\;  H^1(K, S) & \lra  &  H^1(L, S)  & \lga  &  H^1(k(C), S) \\
\quad  \downarrow   &    & \quad \downarrow &    & \quad \downarrow \\
\quad H^1(K, G) & \lra  & \;  H^1(L, G)  & \lga  & \;  H^1(k(C), G) \,
\end{array}  \]
shows that there exists an \'etale algebra $E'/k(C)$
and an Azumaya algebra ${\cal A}'/E'$ such that
$$
E \otimes_K L 
{\buildrel \sim \over \lra}
E' \otimes_{k(C)} L, 
\quad \text{and} \quad
{\cal A}'   \otimes_{E'} (E' \otimes_{k(C)} L)
{\buildrel \sim \over \lra}
 \bigl( {\cal A} \otimes_{E} (E\otimes_K L) \bigr).
$$
Since ${\rm cd}(k(C)) \leq 1$ (see \cite[\S II.3]{serre-gc}),
${\cal A}'/ E'$ is the split Azumaya algebra of rank $n$.
We conclude that
${\cal A} \otimes_E  (E \otimes_K L)/ (E \otimes_K L)$
is the split Azumaya algebra of rank $n$.
This shows that the map
$\Br(E) \rightarrow \Br(E \otimes K_{ab})$ is trivial
for any \'etale algebra $E/K$. Lemma~\ref{lem.brauer} now
tells us that $\cd(K_{ab}) \le 1$. This concludes the proof 
of Theorem~\ref{main}(i).

\smallskip
The proof of part (ii) is exactly the same, except that 
the field extension $L/K$, constructed at the beginning of 
previous paragraph, is now solvable, rather than abelian.
\qed

\begin{remark} \label{rem.hilbert(p)}
A similar argument shows that the conjecture of 
Bogomolov ~\cite[Conjecture 2]{Bog} mentioned in Remark~\ref{rem.main}
is a consequence of the following weaker form of 
Problem~\ref{prob.hilbert}(i) (which is also open):

\smallskip
{\bf Problem~\ref{prob.hilbert}$'$:}
{\em Let $k$ be an algebraically closed field of characteristic zero,
$S$ be a finite group, $K/k$ be a field extension and $p$ be a prime
integer.  Is it true that for every $\alpha \in H^1(K, S)$ there exists
a finite extension $[K' : K]$ of degree prime to $p$ and
an abelian extension $L/K'$ such that $\ed(\alpha_L) \le 1$?}
\end{remark}

\section*{Acknowledgements}  
We are grateful to M. Artin, G. Berhuy, J. Buhler,
J - L. Colliot-Th\'el\`ene, K. Karu, R. Parimala, 
V. L. Popov, J.-P. Serre and B. Youssin for stimulating 
discussions. We also thank M. Florence for pointing out a mistake
in our earlier proof of Proposition~\ref{prop1}. This mistake
has been corrected in the current version.

\end{document}